\documentclass[10pt]{article}
\usepackage{amsmath}
\usepackage{amssymb}
\usepackage{amsthm}
\usepackage{amsfonts}
\usepackage{url}
\theoremstyle{plain}
\newtheorem{theorem}{Theorem}[section]

\newtheorem{lemma}[theorem]{Lemma}
\newtheorem{corollary}[theorem]{Corollary}

\theoremstyle{definition}

\newtheorem{example}{Example}[section]

\theoremstyle{remark}

\newcommand{\cA}{{\mathcal A}}

\newcommand{\NN}{{\mathbb N}}

\newcommand{\cAstar}{\mathcal{A}^*}
\newcommand{\cAplus}{\mathcal{A}^+}
\newcommand{\empt}{\varepsilon}
\newcommand{\cAw}{\mathcal{A}^{\omega}}
\newcommand{\cAinf}{\mathcal{A}^{\infty}}

\newcommand{\rev}{\widetilde}

\newcommand{\vvqed}{\vspace{-1.5cm} \flushright{\qedsymbol}}

\newcommand{\ee}{\'{e}\'{e}}

\newcommand{\ov}{\overline}

\author{Amy Glen\footnotemark[1]}
\title{Conjugates of characteristic Sturmian words \\ generated by morphisms}
\date{December 15, 2003}
\begin{document}
\normalsize
\maketitle 
\footnotetext[1]{E-mail:
\texttt{amy.glen@adelaide.edu.au}; 
URL: \texttt{http://www.maths.adelaide.edu.au/$\sim$aglen}}
\begin{center}
\vspace{-0.7cm} School of Pure Mathematics, University of
Adelaide, South Australia, 5005
\end{center}
\vspace{0.5cm}
\begin{abstract}
This article is concerned with characteristic Sturmian words of
slope $\alpha$ and $1 - \alpha$ (denoted by $c_\alpha$ and $c_{1 -
\alpha}$ respectively), where $\alpha \in (0,1)$ is an irrational
number such that $\alpha = [0;1+d_1,\overline{d_2,\ldots,d_n}]$ with 
$d_n \geq d_1 \geq 1$.
It is known 
that both $c_\alpha$ and $c_{1 - \alpha}$ are fixed points of 
non-trivial (standard) morphisms
$\sigma$ and $\hat{\sigma}$, respectively, if and only if $\alpha$
has a continued fraction expansion as above. Accordingly, such
words $c_\alpha$ and $c_{1 -\alpha}$ are \emph{generated} by the
respective morphisms $\sigma$ and $\hat{\sigma}$. 
For the particular
case when $\alpha = [0;2,\overline{r}]$ ($r\geq1$), we give a decomposition
of each \emph{conjugate} of $c_\alpha$ (and hence $c_{1-\alpha}$)
into \emph{generalized adjoining singular words}, by considering
conjugates of powers of the standard morphism $\sigma$ by which it
is generated. This extends a recent result of Lev\'{e} and S\ee bold
on conjugates of the infinite Fibonacci word. 
\vspace{0.2cm}
\vspace{0.2cm}

\end{abstract}

\section{Introduction}

In recent years, combinatorial properties of finite and infinite
words have become significantly important in fields of physics,
biology, mathematics, and computer science. In particular, the
fascinating family of Sturmian words has been the subject of 
many papers (see \cite{jBpS02stur, gHmM40symb, gM99lynd}, for
example). These words, which represent the simplest family of
quasi-crystals, have numerous applications in various fields of
mathematics, such as symbolic dynamics, the study of continued
fraction expansion, and also in some domains of physics (crystallography) 
and computer science (formal language theory, algorithms on
words, pattern recognition).

  Sturmian words are (aperiodic) infinite words with exactly $n+1$
distinct factors of length $n$, for each $n \in \NN$. Since this
implies that a Sturmian word has exactly two factors of length 1, then any 
such word is over a two-letter alphabet, say $\cA = \{a,b\}$.
There are many characterizations and numerous properties of
Sturmian words. For a comprehensive study of the basic properties
of Sturmian words, and of their transformations by morphisms, see
Berstel and S\ee bold \cite{jBpS02stur}.

  Here, an \emph{infinite word} $x$ over an alphabet $\cA$ is a map 
$x: \NN \longrightarrow
\cA$. For any $i \geq 0$, we set $x_i = x(i)$ and write $x = x_0x_1x_2\cdots$,
where each $x_i \in \cA$. In this paper, we will utilize the following 
characterization of Sturmian words, which was originally proved by Morse and 
Hedlund \cite{gHmM40symb}. An infinite word $s$ over $\cA =
\{a,b\}$ is Sturmian if and only if there exists an irrational
$\alpha \in (0,1)$, and a real number $\rho$, such that $s$ is one
of the following two infinite words: \vspace{-0.4cm}
\[
  s_{\alpha,\rho}, ~s_{\alpha,\rho}^{\prime}: \NN \longrightarrow \cA
\]
defined by
\[
 \begin{matrix}
  &s_{\alpha,\rho}(n) = \begin{cases}
                        a    ~~~~~\mbox{if} ~\lfloor(n+1)\alpha + \rho\rfloor -
                        \lfloor n\alpha + \rho\rfloor = 0, \\
                        b    ~~~~~\mbox{otherwise};
                       \end{cases} \\
  &\qquad \\
  &s_{\alpha,\rho}^\prime(n) = \begin{cases}
                        a    ~~~~~\mbox{if} ~\lceil(n+1)\alpha + \rho\rceil -
                        \lceil n\alpha + \rho\rceil = 0, \\
                        b    ~~~~~\mbox{otherwise}.
                       \end{cases}
 \end{matrix} \qquad (n \geq 0)
\]

The irrational $\alpha$ is called the \emph{slope} of $s$ and
$\rho$ is the \emph{intercept}. If $\rho = 0$, we have
\[
  s_{\alpha,0} = ac_{\alpha} ~\mbox{and} ~s_{\alpha,0}^\prime =
  bc_{\alpha},
\]
where $c_{\alpha}$ is called the \emph{characteristic Sturmian
word} of slope $\alpha$. 

The infinite Fibonacci word $f$ is a special example of a
characteristic Sturmian word of slope \\ $\alpha =
[0;2,\overline{1}] = (3 - \sqrt5)/2$, which is generated by a 
(standard) morphism.
Wen and Wen \cite{zWzW94some} have established a factorization of
the Fibonacci word into \emph{singular words} and, in a similar
fashion, Melan\c{c}on \cite{gM99lynd} has proposed a
generalization of singular words over a two-letter alphabet that
allows for a decomposition of all the characteristic Sturmian
words. More recently,
Lev\'{e} and S\ee bold \cite{fLpS03conj} have obtained a
generalization of Wen and Wen's `singular' decomposition of the Fibonacci
word, by establishing a similar decomposition for each \emph{conjugate}
of this infinite word into what they call \emph{generalized
singular words}. The aim of this current paper is to extend this
latter result to any characteristic Sturmian word of slope $\alpha
= [0;2,\overline{r}]$ (resp. $1-\alpha =[0;1,1,\overline{r}]$), $r \geq 1$,
which we will show is generated by a particular \emph{standard}
morphism $\sigma$ (resp. $\hat{\sigma}$). 

This paper is organized in the following manner. In Section
\ref{S:prelim}, after recalling some combinatorial notions used in the
study of words and morphisms
(\S \ref{SS:words}), we will consider right conjugation of
standard morphisms (\S \ref{SS:conj}). 
We shall then discuss characteristic Sturmian
words $c_\alpha$ and a `singular' decomposition of such words, which
we will later generalize to each conjugate of $c_\alpha$ for particular 
$\alpha$. 
In the section to follow (\S
\ref{S:conj_of_c_alpha}), we describe all irrationals $\alpha \in (0,1)$  
such that $c_\alpha$ is generated by a morphism, and subsequently obtain
generalizations of Lev\'{e} and S\ee bold's \cite{fLpS03conj} results 
(on conjugates of the Fibonacci word) for $c_\alpha$ and $c_{1 - \alpha}$
with $\alpha = [0;2,\overline{r}]$.

\section{Preliminaries} \label{S:prelim}

\subsection{Words and Morphisms} \label{SS:words}

Any of the following terminology that is not further clarified can be
found in either \cite{mL83comb} or \cite{jBpS02stur}, which give more
detailed presentations.

Let $\cA$ be a finite set of symbols that we shall call an
\emph{alphabet}, the elements of which are called \emph{letters}.
A (finite) \emph{word} is an element of the \emph{free monoid} $\cAstar$
generated by $\cA$, in the sense of concatenation. 
The identity $\empt$ of $\cAstar$ is called the \emph{empty word}, and
the \emph{free semi-group}, denoted by $\cAplus$, is defined by 
$\cAplus := \cAstar\setminus\{\empt\}$. We denote by
$\cAw$ the set of all infinite words over $\cA$, and define
$\cAinf := \cAstar \cup \cAw$.

A finite word $w$ is a \emph{factor} of $x \in \cAinf$ if $x = uwv$
for some $u \in \cAstar$ and $v \in \cAinf$. Furthermore, $w$ is
called a \emph{prefix} (resp. \emph{suffix}) of $x$ if $u = \empt$
(resp. $v = \empt$). An infinite word $z \in \cAw$ is called a 
\emph{suffix} of $x \in \cAw$ if there is a word $w \in \cAstar$ 
such that $x = wz$. The \emph{length} $|w|$ of
a finite word $w$ is defined to be the number of letters it contains. 
(Note that $|\empt|=0$.)

The \emph{inverse} of $w \in \cAstar$, written $w^{-1}$, is
defined by $ww^{-1} = w^{-1}w = \empt$. It must be emphasized that
this is merely notation, i.e. for $u, v, w \in \cAstar$, the words
$u^{-1}w$ and $wv^{-1}$ are defined only if $u$ (resp. $v$) is a
prefix (resp. suffix) of $w$. Also note that if $w = uv$ then
$wv^{-1} = u$ and $u^{-1}w = v$, and if $x = wx^\prime$, where $w
\in \cAstar$ and $x^\prime \in \cAw$, then $w^{-1}x = x^\prime$.

Two words $w$, $z \in \cAstar$ are said to be \emph{conjugate} if
there exist words $u$, $v$ such that $w = uv$ and $z = vu$. If
$|u| = k$, then $z$ is called the $k$-th conjugate of $w$. This
notion extends to infinite words as follows. For $k \in \NN$, the
$k$-th \emph{conjugate} of an infinite word $x$ over $\cA$ is
the infinite word $x^\prime$ such that $x = ux^ \prime$, where $u
\in \cAstar$ and $|u| = k$.

A \emph{morphism on} $\cA$ is
a map $\psi: \cAstar \longrightarrow \cAstar$ such that $\psi(uv)
= \psi(u)\psi(v)$ for all $u, v \in \cAstar$.  It is uniquely
determined by its image on the alphabet $\cA$. If $\psi(c) = cw$,
for some letter $c \in \cA$ and some $w \in \cA^{+}$, then $\psi$ is
said to be \emph{prolongable} on $c$. In this case, the word
$\psi^{n}(c)$ is a proper prefix of the word $\psi^{n+1}(c)$ for
each $n \in \NN$, and the sequence $(\psi^n(c))_{n\geq 0}$
converges to a unique infinite word
\[
  x = \underset{n \rightarrow\infty}{\lim}\psi^n(c) = \psi^{\omega}(c).
\]
An infinite word $x$ is \emph{generated by a morphism} if $x =
\psi^\omega(c)$ for some letter $c$ and some morphism $\psi$.

In what follows, it is assumed that all words are over the
two-letter alphabet $\cA = \{a,b\}$.

\subsection{Conjugation of Standard Morphisms} \label{SS:conj}

Define on $\cA$ the following three morphisms
\[
  E: \begin{matrix}
      &a &\mapsto &b& \\
      &b &\mapsto &a&
     \end{matrix}, \qquad \varphi: \begin{matrix}
                                   &a &\mapsto &ab& \\
                                   &b &\mapsto &a& 
                                   \end{matrix}, \qquad 
  \rev{\varphi}: \begin{matrix}
                  &a &\mapsto &ba& \\
                  &b &\mapsto &a~&
                 \end{matrix}.
\]
A morphism $\psi$ is \emph{Sturmian} if and only if $\psi \in
\{E,\varphi,\rev{\varphi }\}^{*}$, i.e. if and only if it is a
composition of $E$, $\varphi$, and $\rev{\varphi}$ in any number
and order (see \cite{fMpS93morp}). Furthermore, 
a morphism $\psi$ is \emph{standard} if and only if $\psi \in
\{E,\varphi\} ^{*}$ (see \cite{aD97stan}). Note that a morphism is 
\emph{non-trivial} if it is neither $E$ nor $Id_{\cA}$ (the identity morphism).

Suppose $\psi$ and $\xi$ are morphisms on $\cA$. If there exists a word 
$u$ such that
\[ 
  \psi(w)u = u\xi(w) \quad \mbox{for all words $w \in \cAstar$}, 
\]
then $\xi$ is called the $|u|$-th \emph{right conjugate} of $\psi$, 
denoted $\psi_{|u|}$.

It has been shown by S\ee bold \cite{pS98onth} that
the number of distinct right conjugates of a standard morphism
$\psi$ is $|\psi(ab)| - 1$; namely, the morphisms $\psi_{0}$ to
$\psi_{|\psi(ab)|-2}$. The following useful lemma is proved in 
\cite{fLpS03conj}.

\begin{lemma} \label{L:2.2}
Suppose the infinite word $x$ is generated by the standard morphism $\psi$.
Let $k \in \NN$ with \\ $0 \leq k \leq |\psi(ab)| - 2$, and let $v$
denote the prefix of $x$ of length $k$. Then
$\psi_k$ is such that $\psi_k(x) = v^{-1}x$. 
\qed
\end{lemma}

Thus, if $\psi$ is a
standard morphism that generates an infinite word $x$, one deduces
from Lemma \ref{L:2.2} that the result of applying $\psi_k$ to $x$
simply consists of deleting the first $k$ letters of $x$, i.e. $\psi_k(x)$
is the $k$-th conjugate of $x$.

\subsection{Characteristic Sturmian Words $c_\alpha$ and Singular Words}
\label{SS:Sturmian}

Note that every irrational $\alpha \in (0,1)$ has a unique
continued fraction expansion
\[ 
  \alpha = [0;a_1,a_2,a_3,\ldots] = \cfrac{1}{a_1+
                                \cfrac{1}{a_2 +
                                \cfrac{1}{a_3 + \cdots,
                                 }}}
\]
where each $a_i$ is a positive integer. If the sequence
$(a_i)_{i\geq1}$ is eventually periodic, with $a_i = a_{i+m}$ for
all $i \geq n$, we use the notation
$\alpha = [0;a_1,a_2,\ldots,a_{n-1},\ov{a_{n},a_{n+1},\ldots,a_{n+m-1}}]$.
The $n$-th \emph{convergent} of $\alpha$ is defined by
\[
  \frac{p_n}{q_n} = [0;a_1,a_2,\ldots,a_n] \quad \mbox{for all} ~n\geq 1,
\]
where the sequences $(p_n)_{n\geq0}$ and $(q_n)_{n\geq0}$ are
given by
\[
\begin{matrix}
 &p_{0} = 0, &p_{1} = 1, &p_n = a_np_{n-1} + p_{n-2}, ~~&n\geq 2; \\
 &q_{0} = 1, &q_{1} = a_1, &q_n = a_nq_{n-1} + q_{n-2}, ~~&n\geq 2.
\end{matrix}
\]

Suppose $\alpha = [0;1+d_1,d_2,d_3, \ldots]$, with $d_1 \geq 0$
and all other $d_n > 0$. To the \emph{directive sequence}
$(d_1,d_2,d_3,\ldots)$, we associate a sequence $(s_n)_{n \geq
-1}$ of words defined by
\[
  s_{-1} = b, ~s_{0} = a, ~s_{n} = s_{n-1}^{d_{n}}s_{n-2}; \quad n \geq 1.
\]
Such a sequence of words is called a \emph{standard sequence}, and
we have
\[
  |s_n| = q_n \quad \mbox{for all} ~n\geq0.
\]
Note that $ab$ is a suffix of $s_{2n-1}$ and $ba$ is a suffix of
$s_{2n}$, for all $n \geq 1$.

Standard sequences are related to characteristic Sturmian words in
the following way. Observe that, for any $n\geq0$, $s_n$ is a
prefix of $s_{n+1}$, which gives obvious meaning to $\underset{n
\rightarrow \infty}{\mbox{lim}}s_n$ as an infinite word. In
fact, one can prove \cite{aFmMuT78dete,tB93desc} that each $s_n$
is a prefix of $c_\alpha$, and we have
\[ 
c_{\alpha} = \underset{n \rightarrow \infty}{\mbox{lim}}s_n.
\]

\subsubsection{Singular Decomposition of $c_\alpha$}

Melan\c{c}on \cite{gM99lynd} (also see \cite{wCzW03some}) has
proposed a generalization of Wen and Wen's \cite{zWzW94some} 
singular factors of the Fibonacci word to the case of any characteristic
Sturmian word $c_\alpha$ and, in doing so, has established a decomposition
of $c_\alpha$ into \emph{adjoining singular words}, as shown
below.

For $c_\alpha$ such that $\alpha = [0;1+d_1,d_2,d_3, \ldots]$ with
$d_1 \geq 1$, Melan\c{c}on \cite{gM99lynd} introduced the singular
words $w_n$ of $c_\alpha$ defined by
\[
  w_n = \begin{cases}
         as_{n}b^{-1} \qquad \mbox{if $n$ is odd}, \\
         bs_na^{-1} \qquad \mbox{otherwise},
        \end{cases}
\] for $n \geq 1$, with $w_{-2}= \empt, w_{-1} =
a, w_{0} = b$. Furthermore, the following words $v_n$ are also
defined in \cite{gM99lynd}. For all $n \geq -1$,
\[
  v_n = \begin{cases} as_{n+1}^{d_{n+2} - 1}
                          s_n b^{-1} ~~~\mbox{if $n$ is odd}, \\
                          bs_{n+1}^{d_{n+2} - 1}
                          s_n a^{-1} ~~~\mbox{otherwise}.
                          \end{cases}
\]
Clearly, the word $v_n$ differs from $w_{n+2}$ by a factor
$s_{n+1}$, and it can be proved that all $v_n$ and $w_n$ are palindromes
(i.e. words that read the same backwards as forwards).
We shall call $v_n$ the
$n$-th \emph{adjoining singular word} of $c_{\alpha}$,  
and set $v_{-2} = \empt$.
In terms of the singular and adjoining singular words of
$c_\alpha$, the following generalization of Wen and Wen's
\cite{zWzW94some} singular decomposition of the Fibonacci word 
has been established.
\begin{theorem} \emph{\cite{gM99lynd}} \label{T:melancon1}
~$c_\alpha = \prod_{j=-1}^\infty(v_{2j}w_{2j+1})^{d_{2j+3}}=
        \prod_{j=-1}^\infty v_{j}$.\qed
\end{theorem}

In the next section, for the case $\alpha =
[0;2,\ov{r}]$, we will generalize this factorization of $c_\alpha$
(and hence $c_{1-\alpha}$), by showing that, for each prefix $v$
of $c_\alpha$, $v^{-1}c_\alpha$ can be decomposed into
\emph{generalized adjoining singular words}. Such a result has
already been established (by Lev\'{e} and S\ee bold
\cite{fLpS03conj}) for the case of the Fibonacci word $f =
c_{(3-\sqrt{5})/2}$, where $\frac{3-\sqrt{5}}{2} =
[0;2,\ov{1}]$. 

\section{Decompositions of Conjugates of $c_\alpha$} \label{S:conj_of_c_alpha}

\subsection{Characteristic Sturmian Words Generated by Morphisms}

Here, we describe all irrationals $\alpha \in (0,1)$ such that the
characteristic Sturmian word $c_\alpha$ is generated by a
morphism.
In order to do this, we must first define a special set
of irrational numbers. A \emph{Sturm number} (see
\cite{jBpS02stur}) is an irrational number $\alpha \in (0,1)$ that
has a continued fraction expansion of one of the following types:
\begin{itemize}
\item[(i)] ~$\alpha = [0;1 + d_1,\overline{d_2,\ldots,d_n}] <
\frac{1}{2}$ with ~$d_n \geq d_1 \geq 1$;
\item[(ii)] ~$\alpha = [0;1,d_1,\overline{d_2,\ldots,d_n}] >
\frac{1}{2}$ with ~$d_n \geq d_1$.
\end{itemize}
Observe that if $\alpha = [0;1+d_1,\overline{d_2,\ldots,d_n}]$
with $d_n \geq d_1 \geq 1$, then
\[
  1 - \alpha = \frac{1}{1 + \alpha/(1-\alpha)} = [0;1,d_1,\overline{d_2,
                        \ldots,d_n}].
\]
Hence, $\alpha$ has an expansion of type (i) if and only if $1 -
\alpha$ has an expansion of type (ii). Accordingly, $\alpha$ is
a Sturm number if and only if $1-\alpha$ is a Sturm number.

In what follows, we will always assume (unless otherwise stated) that
$\alpha$ is a Sturm number of type (i). Also, we shall denote the 
standard sequence of
$c_{\alpha}$ (resp. $c_{1-\alpha}$) by $(s_n)_{n\geq-1}$ (resp.
$(\hat{s}_n)_{n\geq-1}$). Clearly, we have $\hat{s}_{1} = \hat{s}_{-1} = b$ 
since $1 - \alpha = [0;1,d_1,\ov{d_2,\ldots,d_n}]$. Consequently,
$c_{1-\alpha}$ is obtained from $c_{\alpha}$ by exchanging all
letters $a$ and $b$ in $c_\alpha$, i.e. $c_{1-\alpha} =
E(c_{\alpha})$. Indeed, it is easily checked that
\[
  \hat{s}_n = E(s_{n-1}) \quad \mbox{for all} ~n\geq0.
\]
Hence, 
\[
  E(c_{\alpha}) = E\left(\underset{n \rightarrow \infty}{\mbox{lim}}s_{n}
                     \right) \\
       = \underset{n \rightarrow \infty}{\mbox{lim}}E(s_{n}) 
       = \underset{n \rightarrow \infty}{\mbox{lim}}\hat{s}_{n+1} 
       = c_{1-\alpha}.
\]
Therefore, we can restrict our attention to
characteristic Sturmian words $c_\alpha$ such that $\alpha$ is a 
Sturm number of type (i). Later, an analogue of the 
main result of this paper (Theorem \ref{T:12.06.03(1)}) will be deduced for
$c_{1-\alpha}$. 

We say that a morphism $\psi$ \emph{fixes} an infinite word $x$ if
$\psi(x) = x$, in which case $x$ is called a \emph{fixed point} of
$\psi$. The following result describes all irrationals $\alpha \in (0,1)$ such
that $c_\alpha$ is a fixed point of a non-trivial morphism.

\begin{theorem} \emph{\cite{dCwMaPpS93subs, tKaV96subs, jBpS93acha}}
\label{T:tKaV96subs} Let $\alpha \in (0,1)$ be irrational. Then
$c_\alpha$ is a fixed point of a non-trivial morphism $\sigma$ if
and only if $\alpha$ is a Sturm number. In particular, if $\alpha
= [0;1 + d_1,\overline{d_2,\ldots,d_n}]$ with ~$d_n
\geq d_1 \geq 1$, then $c_\alpha$ is the fixed point of any power
of the morphism
\[
  \sigma: \begin{matrix}
           &a& &\mapsto &s_{n-1}~~~~~~~& \\
           &b& &\mapsto &s_{n-1}^{d_n - d_1}s_{n-2}&
          \end{matrix}.
\]
Further, $c_{1-\alpha}$ is a fixed point of any power of the
morphism
\[
  \hat{\sigma}: \begin{matrix}
                 &a& &\mapsto &\hat{s}_{n}^{d_n - d_1}\hat{s}_{n-1}& \\
                 &b& &\mapsto &\hat{s}_{n}~~~~~~~~~~&
                \end{matrix}.
\]\qed
\end{theorem}

Note that, for any $m \in \NN$, both $\sigma^m$ and $\hat{\sigma}^m$ are 
standard morphisms. In fact, it was shown by Crisp et al.
\cite{dCwMaPpS93subs} that \vspace{-0.4cm}
\begin{align*} 
  \sigma &= (\varphi E)^{d_1}E(\varphi E)^{d_2}E\cdots
  (\varphi E)^{d_{n-1}}E(\varphi E)^{d_n - d_1}; ~~\mbox{and} \\
  \hat{\sigma} &= E(\varphi E)^{d_1}E(\varphi E)^{d_2}E\cdots
  (\varphi E)^{d_{n-1}}E(\varphi E)^{d_n - d_1}E.
\end{align*}
Observe that $\hat{\sigma} = E\sigma E$.

Now, S\ee bold \cite{pS98onth} proved that a standard morphism
$\psi$ generates an infinite (characteristic Sturmian) word if and
only if
\[
  \psi \in \{\varphi,E\varphi,\varphi E,E\varphi E\}^{+}\backslash
(\{E\varphi\}^+ \cup \{\varphi E\}^+).
\]
Here, we prove that a characteristic Sturmian word $c_{\gamma}$ is
generated by a (standard) morphism if and only if $\gamma$ is a
Sturm number. Specifically, we prove that $c_\alpha =
\sigma^\omega(a)$ and $c_{1-\alpha} = \hat{\sigma}^\omega(b) =
(E\sigma E)^\omega(b)$. These are direct results of the following
lemma and corollary.

\begin{lemma} \label{L:10.06.03(1)}
For any $k \in \NN$, $\sigma(s_k) = s_{k + (n-1)}.$
Consequently, if $k \in \NN$ is fixed, then
\[
  \sigma^m(s_k) = s_{k + m(n-1)} \quad \mbox{for all} ~m \geq 0.
\]
\end{lemma}
\begin{proof} Mathematical induction.
\end{proof}

\begin{corollary} \label{Cor:10.06.03(2)}
For any integer $m \geq 1$,
 $\sigma^m(a) = s_{m(n-1)} 
 ~~\mbox{and} ~~\sigma^m(b) = s_{m(n-1)}^{d_n - d_1}s_{m(n-1)-1}$. 
\qed
\end{corollary}
\vspace{0.2cm}
As an immediate consequence of the above corollary, we have the
following result.

\begin{corollary} \label{Cor:07.07.03(1)} Let
$\alpha = [0;1+d_1,\overline{d_2, \cdots,d_n}]$ with $d_n \geq d_1
\geq 1$. Then
\begin{itemize}
\item[\emph{(i)}] 
~$c_\alpha = \underset{m \rightarrow \infty}{\mbox{\emph{lim}}}
               \sigma^m(a) = \sigma^\omega(a);$
\item[\emph{(ii)}] ~$c_{1-\alpha} = \underset{m \rightarrow \infty}{\mbox
{\emph{lim}}} \hat{\sigma}^m(b) = \hat{\sigma}^\omega(b)$,
~\mbox{where $\hat{\sigma} = E\sigma E$}.
\end{itemize}
\end{corollary}
\begin{proof}
The fact that $c_\alpha = \underset{m \rightarrow
\infty}{\mbox{lim}} \sigma^m(a)$ follows from Corollary
\ref{Cor:10.06.03(2)} since $\sigma^m(a) = s_{m(n-1)}$, for any
integer $m \geq 1$. Moreover, we know that $c_{1-\alpha} =
E(c_{\alpha})$, so that (ii) is obtained by realizing
\[
  c_{1-\alpha} = E\left(\underset{m \rightarrow \infty}{\mbox{lim}}
                  \sigma^m(a)\right)
               = \underset{m \rightarrow \infty}{\mbox{lim}}E\sigma^mE(b)
               = \underset{m \rightarrow \infty}{\mbox{lim}}(E\sigma E)^m(b).
\]
\end{proof}

In light of Corollary 
\ref{Cor:07.07.03(1)}, one has that if $\gamma$ is a Sturm number,
then $c_\gamma$ is generated by a (standard) morphism. 
The converse is trivially true in view of Theorem \ref{T:tKaV96subs}.

By considering Melan\c{c}on's factorization of $c_\alpha$ into 
adjoining singular words (Theorem \ref{T:melancon1}), 
we will now extend Lev\'{e} and S\ee bold's 
result (Theorem 4.6 in \cite{fLpS03conj}) to the case 
$\alpha = [0;2,\ov{r}]$.

\subsection{The case $\alpha = [0;2,\ov{r}]$} \label{SS:alpha*}

Now, if $\alpha = [0;2,\overline{r}]$, then for each $m \in \NN$,
\[
   v_m = \begin{cases}
           as_{m+1}^{r-1}s_mb^{-1} ~~~\mbox{if $m$ is odd}, \\
           bs_{m+1}^{r-1}s_ma^{-1} ~~~\mbox{otherwise},
        \end{cases}
\] and $v_{-1} = as_0^0s_{-1}b^{-1} = a$.
Observe that, for any integer $m \geq 0$,
$|\sigma^m(ab)| = |\sigma(s_1)| = |s_{m+1}| = q_{m+1}$.
Furthermore, using Corollary \ref{Cor:10.06.03(2)}, it is 
easily checked that, for any $m \geq -1$,
\[
  v_m = \begin{cases}
           a\sigma^{m+1}(b)b^{-1} ~~~\mbox{if $m$ is odd}, \\
           b\sigma^{m+1}(b)a^{-1} ~~~\mbox{otherwise},
        \end{cases}
\]
since $\sigma^{m+1}(b) = s_{m+1}^{r-1}s_{m}$, for all $m \in \NN$.
Also note that
$|v_m| = q_{m+2} - q_{m+1}, ~\mbox{for all $m \geq -1$}$.

Whence, we have the following special case of Lemma
\ref{L:2.2}.

\begin{lemma} \label{L:11.06.03(1)}
Suppose $\alpha = [0;2,\overline{r}]$ and let $k, m \in \NN$ be
such that $0 \leq k \leq q_{m+1} - 2$. If $v$ denotes the prefix
of length $k$ of $c_\alpha (= (\sigma^m)^\omega(a))$, then
$(\sigma^{m})_k(c_\alpha) = v^{-1}c_\alpha$. \qed
\end{lemma}
\vspace{0.2cm}
The next lemma shows how to remove a prefix from the `singular' 
decomposition of $c_\alpha$ \\ 
(\emph{cf.} Proposition 4.5 in \cite{fLpS03conj}).

\begin{lemma} \label{T:11.06.03(2)}
Suppose $\alpha = [0;2,\overline{r}]$ and let $k, m \in \NN$  be
such that $k = q_{m+1} - p$ with ~$2 \leq p \leq q_{m+1} - q_m +
1$. Then
\[
  (\sigma^m)_k(c_\alpha) = u^{-1}v_{m-1}\prod_{j=m}^\infty v_{j},
\]
where $u$ is the prefix of $v_{m-1}$ of length $|u| = q_{m+1} -
q_{m} + 1 - p$.
\end{lemma}
\begin{proof}
We have $q_m - 1 \leq k \leq q_{m+1} - 2$. Thus,
if $k = 0$, then we must have $m  = 0$. Now, $k = q_1 - p = 0$
implies $p = q_1 = 2$, and therefore, $|u| = q_1 - q_0 + 1 - 2 =
0$. Further, from Theorem \ref{T:melancon1}, we have
\[
  (\sigma^0)_0(c_\alpha) = c_\alpha = \prod_{j=-1}^{\infty}v_{j},
\]
so the result holds for $k = 0$.

Suppose $k \geq 1$, then $m \geq 1$ and, in this case, observe
that
\[
  q_m - 1 = \sum_{j=-1}^{m-2}(q_{j+2} - q_{j+1}).
\]
From Lemma \ref{L:11.06.03(1)}, we know that
$(\sigma^{m})_k(c_\alpha)$ is the word obtained from $c_\alpha$ by
removing its prefix of length $k$, i.e. of length at least
$\sum_{j=-1}^{m-2} (q_{j+2} - q_{j+1})$. Whence, since 
$|v_j| = q_{j+2} - q_{j+1}$ for any $j \geq -1$, then
$(\sigma^m)_k(c_\alpha)$ is obtained from $c_\alpha = \prod_
{j=-1}^\infty v_j$  by first removing the prefix
$v_{-1}v_0v_1\cdots v_{m-2}$. Then, from the remaining
infinite word $v_{m-1}\prod_{j=m}^{\infty}v_j$, we remove the
prefix $u$ of $v_{m-1}$ of length
\[
  |u| = k - (q_{m} - 1) = q_{m+1} - p - q_m + 1 = q_{m+1} - q_m + 1 - p.
\]
\end{proof}

\begin{example} \label{ex:12.06.03}
Take $\alpha = [0;2,\overline{3}] = (\sqrt{13} - 1)/6$, so that
\[
  c_\alpha = \underline{a}~\underline{babab}~\underline{aabababaabababaa}
             ~babababaabababaabababaabababab\cdots
\]
Note that $v_{-1} = a, v_0 = babab, v_1 = aabababaabababaa$
since $s_1 = ab$ and $s_2 = abababa$. Hence, by the preceding
lemma,
\begin{alignat*}{3}
c_\alpha &= &~(\sigma^0)_0(c_\alpha) &= v_{-1}v_0v_1v_2v_3\cdots \\
a^{-1}c_\alpha &= &~(\sigma^1)_1(c_\alpha) &= v_0v_1v_2v_3\cdots \\
(ab)^{-1}c_\alpha &= &~(\sigma^1)_2(c_\alpha) &=
b^{-1}v_0v_1v_2v_3
                      \cdots \\
(aba)^{-1}c_\alpha &= &~(\sigma^1)_3(c_\alpha) &=
(ba)^{-1}v_0v_1v_2
                      v_3\cdots \\
(abab)^{-1}c_\alpha &= &~(\sigma^1)_4(c_\alpha) &=
(bab)^{-1}v_0v_1
                        v_2v_3\cdots \\
(ababa)^{-1}c_\alpha &= &~(\sigma^1)_5(c_\alpha) &=
(baba)^{-1}v_0v_1
                        v_2v_3\cdots \\
(ababab)^{-1}c_\alpha &= &~(\sigma^2)_6(c_\alpha) &=
v_1v_2v_3\cdots,
                ~\mbox{etc}.
\end{alignat*} 
\end{example}
\vspace{0.2cm}
For any $n \geq -1$, set $(\sigma^{n+1})_{-1}(b) = v_n$.

\begin{theorem} \label{T:12.06.03(1)}
Let $k, m \in \NN$ be such that $k = q_{m+1} - p$ with ~$2 \leq p
\leq q_{m+1} - q_{m} + 1$. Then
\[
  (\sigma^m)_k(c_\alpha) = \prod_{j=m-1}^\infty (\sigma^{j+1})_
{q_{m+1} - q_m - p}(b).
\]
\end{theorem}
\begin{proof}
It follows immediately from Lemma \ref{T:11.06.03(2)}  that
\[
  (\sigma^{m})_k(c_\alpha) = u^{-1}v_{m-1}\prod_{j=m}^{\infty}v_j
               = u^{-1}\prod_{j=m-1}^{\infty}v_j,
\]
where $u$ is the prefix of $v_{m-1}$ such that $|u| = q_{m+1} -
q_m + 1 - p$. Thus, if $u = \empt$, then
\[
  (\sigma^{m})_k(c_\alpha) = \prod_{j=m-1}^{\infty}v_j
                           = \prod_{j=m-1}^{\infty}(\sigma^{j+1})_{-1}(b),
\] and so the result holds since $|u| = 0 = q_{m+1} - q_m
+ 1 - p$ implies $q_{m+1} - q_m - p = -1$.

If $k = 0$, then $m = 0$, and hence, $p = 2$ so that $|u| =
q_{1} - q_{0} + 1 - 2 = 2 - 1 - 1 = 0$ (i.e. $u =
\empt$). So the result holds for $m=0$, and we therefore take $k
\geq 1$, so that $m \geq 1$.

Observe that, by definition of the adjoining singular word $v_m$, 
there exist letters $x,y \in \cA$ $(x \ne y)$ such
that, for any integer $p \geq -1$,
\[
  v_{m+p} = \begin{cases}
         x\sigma^{m+p+1}(b)y^{-1} ~~~\mbox{if $p$ is odd}, \\
         y\sigma^{m+p+1}(b)x^{-1} ~~~\mbox{otherwise}.
        \end{cases}
\]
Hence, \vspace{-0.4cm}
\begin{align*}
  \prod_{j=m-1}^{\infty}v_j& = (x\sigma^{m}(b)y^{-1})
                             (y\sigma^{m+1}(b)x^{-1})
                             (x\sigma^{m+2}(b)y^{-1})\cdots \\
                          &= x\sigma^{m}(b)\sigma^{m+1}(b)
                             \sigma^{m+2}(b)\cdots \\
                          &= x\prod_{j=m-1}^{\infty}\sigma^{j+1}(b),
\end{align*}
and therefore, $(\sigma^{m})_k(c_\alpha) = u^{-1}x\prod_{j=m-1}^
{\infty}\sigma^{j+1}(b)$. If $u \ne \empt$, then there exists
a word $\hat{u}$ such that $u^{-1}x = {\hat{u}}^{-1}$, which implies that
$\hat{u} = x^{-1}u$ with $|\hat{u}| = q_{m+1} - q_{m} - p$.

For any $n \in \NN$, $\sigma^{n+1}(b) = s_{n+1}^{r-1}s_n$  is a
prefix of $\sigma^{n+2}(b) = s_{n+2}^{r-1}s_{n+1} =
(s_{n+1}^{r-1}s_n)^{r-1}s_{n+1}$, and $|\sigma^{n+1}(b)| = |v_n| =
q_{n+2} - q_{n+1}$. Whence, for any integer $r \geq m - 1$, we
have $\sigma^{r+1}(b) = \hat{u}u_{r+1}$, for some $u_{r+1} \in
\cAstar$ with $|u_{r+1}| \geq p$. Indeed, for $r \geq m - 1$, we have
\begin{align*}
  |u_{r+1}| = |\sigma^{r+1}(b)| - |\hat{u}| 
            &= q_{r+2} - q_{r+1} - q_{m+1} + q_{m} + p \\
            &= (q_{r+2} - q_{m+1}) - (q_{r+1} - q_m) + p \geq p.
\end{align*}

Consequently, by definition of right conjugation of morphisms,
\[
  (\sigma^{r+1})_{q_{m+1} - q_{m} - p}(b) = u_{r+1}\hat{u}.
\]

From the above observations, we therefore find
\begin{align*}
 (\sigma^{m})_k(c_\alpha)&= u^{-1}x\prod_{j=m-1}^
                              {\infty}\sigma^{j+1}(b) \\
                          &= {\hat{u}}^{-1}\prod_{j=m-1}^
                              {\infty}\sigma^{j+1}(b) \\
                          &= {\hat{u}}^{-1}\sigma^m(b)\prod_{j=m}^
                             {\infty}\sigma^{j+1}(b) \\
                          &= {\hat{u}}^{-1}\hat{u} u_m\prod_{j=m}^
                              {\infty}\sigma^{j+1}(b) \\
                          &= u_m\prod_{j=m}^
                              {\infty}\sigma^{j+1}(b) \\
                          &= u_m\hat{u}{\hat{u}}^{-1}\prod_{j=m}^
                              {\infty}\sigma^{j+1}(b) \\
                          &= (\sigma^m)_{q_{m+1} - q_{m} - p}(b)
                             {\hat{u}}^{-1}\prod_{j=m}^
                              {\infty}\sigma^{j+1}(b) \\
                          &= \cdots \\
                          &= \prod_{j=m-1}^
                              {\infty}(\sigma^{j+1})_{q_{m+1} - q_{m}
                             - p}(b).
\end{align*}
\end{proof}

So if $v$ is a prefix of $c_\alpha$, then $v^{-1}c_\alpha$ can be
obtained by concatenating all the words
$[(\sigma^{j+1})_i(b)]_{j\geq l}$, where $i$ and $l$ are integers
depending only on $|v|$. The characteristic Sturmian word
$c_\alpha$ is the special case when $i = l = -1$, so that
\[
  c_\alpha = \prod_{j=-1}^\infty(\sigma^{j+1})_{-1}(b),
\]
where Melan\c{c}on's adjoining singular words, $v_j$, are all the
words $(\sigma^{j+1})_{-1}(b), j \geq -1$.

Recall that the Fibonacci word $f$ is the characteristic Sturmian word
$c_\alpha$ such that $\alpha = [0;2,\overline{1}] = (3 -
\sqrt{5})/2$. In this case, one has $\sigma = \varphi$, i.e.
\[
  f = \varphi^\omega(a) = abaababaabaababaababaabaababaabaab\cdots
\]
If we set $\varphi^{-1}(a) = b$, then for
any integer $n \geq -1$, $\varphi^n(a) = \varphi^{n+1}(b)$ with 
$|\varphi^{n}(a)| = F_n = |\varphi^{n+1}(b)|$, 
where $F_n$ is the $n$-th \emph{Fibonacci
number} defined by
\[
  F_{-1} = F_0 = 1, ~F_n = F_{n-1} + F_{n-2}; \quad n \geq 1.
\]
Note that $q_m = F_m$ for every $m \in \NN$, and
\[
  F_{m+1} - F_m = F_{m+1} - F_m = F_{m-1} \quad \mbox{for all} ~m \in \NN.
\]
Hence, it is deduced from Theorem \ref{T:12.06.03(1)} that if $k,
m \in \NN$ are such that $k = F_{m+1} - p$ with \\ $2 \leq p \leq
F_{m-1} + 1$, then
\[
  (\varphi^m)_k(f) = \prod_{j=m-1}^\infty (\varphi^{j})_
{F_{m-1} - p}(a),
\]
which is Lev\'{e} and S\ee bold's result (Theorem 4.6 in \cite{fLpS03conj}).

As an example, we list some
decompositions of conjugates of the characteristic Sturmian word
$c_\alpha$ for $\alpha = [0;\overline{2}] = \sqrt{2} - 1$.

\begin{alignat*}{14}
&\qquad  &\qquad v_{-1}& \qquad &\quad &v_0 &\quad &v_1 &\quad
&v_2
&\quad &v_3 &\quad \cdots& \\
&m = 0 &\qquad c_\alpha = (\sigma^{0})_{-1}(b)& &\quad
       &(\sigma^{1})_{-1}(b) &\quad
       &(\sigma^{2})_{-1}(b) &\quad &(\sigma^{3})_{-1}(b) &\quad
       &(\sigma^{4})_{-1}(b) &\quad \cdots& \\
&m = 1 &\qquad a^{-1}c_\alpha =& &\quad &(\sigma^{1})_{-1}(b)
&\quad
         &(\sigma^{2})_{-1}(b) &\quad &(\sigma^{3})_{-1}(b) &\quad
         &(\sigma^{4})_{-1}(b) &\quad \cdots& \\
&\qquad &\qquad (ab)^{-1}c_\alpha =& &\quad &(\sigma^{1})_{0}(b)
&\quad
         &(\sigma^{2})_{0}(b) &\quad &(\sigma^{3})_{0}(b) &\quad
         &(\sigma^{4})_{0}(b) &\quad \cdots& \\
&\qquad &\qquad (aba)^{-1}c_\alpha =& &\quad &(\sigma^{1})_{1}(b)
&\quad
         &(\sigma^{2})_{1}(b) &\quad &(\sigma^{3})_{1}(b) &\quad
         &(\sigma^{4})_{1}(b) &\quad \cdots& \\
&m = 2 &\qquad (abab)^{-1}c_\alpha =& &\quad &\qquad &\quad
         &(\sigma^{2})_{-1}(b) &\quad &(\sigma^{3})_{-1}(b) &\quad
         &(\sigma^{4})_{-1}(b) &\quad \cdots& \\
&\qquad &\qquad (ababa)^{-1}c_\alpha =& &\quad &\qquad &\quad
         &(\sigma^{2})_{0}(b) &\quad &(\sigma^{3})_{0}(b) &\quad
         &(\sigma^{4})_{0}(b) &\quad \cdots& \\
&\qquad &\qquad (ababaa)^{-1}c_\alpha =& &\quad &\qquad &\quad
         &(\sigma^{2})_{1}(b) &\quad &(\sigma^{3})_{1}(b) &\quad
         &(\sigma^{4})_{1}(b) &\quad \cdots& \\
&\qquad &\qquad (ababaab)^{-1}c_\alpha =& &\quad &\qquad &\quad
         &(\sigma^{2})_{2}(b) &\quad &(\sigma^{3})_{2}(b) &\quad
         &(\sigma^{4})_{2}(b) &\quad \cdots& \\
&\qquad &\qquad (ababaaba)^{-1}c_\alpha =& &\quad &\qquad &\quad
         &(\sigma^{2})_{3}(b) &\quad &(\sigma^{3})_{3}(b) &\quad
         &(\sigma^{4})_{3}(b) &\quad \cdots& \\
&\qquad &\qquad (ababaabab)^{-1}c_\alpha =& &\quad &\qquad &\quad
         &(\sigma^{2})_{4}(b) &\quad &(\sigma^{3})_{4}(b) &\quad
         &(\sigma^{4})_{4}(b) &\quad \cdots& \\
&\qquad &\qquad (ababaababa)^{-1}c_\alpha =& &\quad &\qquad &\quad
         &(\sigma^{2})_{5}(b) &\quad &(\sigma^{3})_{5}(b) &\quad
         &(\sigma^{4})_{5}(b) &\quad \cdots& \\
&m = 3 &\qquad (ababaababaa)^{-1}c_\alpha =& &\quad &\qquad &\quad
         &\qquad &\quad &(\sigma^{3})_{-1}(b) &\quad
         &(\sigma^{4})_{-1}(b) &\quad \cdots&
\end{alignat*}

\subsection{The case $1 - \alpha = [0;1,1,\ov{r}]$} \label{SS:1-alpha*}

Now, if $\alpha = [0;2,\ov{r}]$, then 
$1 - \alpha = [0;1,1,\ov{r}]$. By observing that $c_{1-\alpha} = 
\hat{\sigma}^\omega(b)$, where $\hat{\sigma} = E\sigma E$, it is
clear that the following theorem is an immediate consequence of
Theorem \ref{T:12.06.03(1)} and the lemma below.

\begin{lemma} \label{L:04.07.03} For any standard morphism $\psi$,
\[
  (E\psi E)_{k} = E\psi_{k}E; \quad 0 \leq k \leq |\psi(ab)| - 2.
\]
\end{lemma}
\begin{proof}
Let $w \in \cAstar$, with $|w| = k$, be such that $\psi(z)w = w
\psi_{k}(z)$ for all $z \in \cA$. Then, for some $z \in \cA$,
\[
  E\psi E(z)E(w) = E(\psi E(z)w) = E(\psi(y)w), \quad \mbox{for some
$y \in \cA$, $y \ne z$}.
\]
Therefore, for $z \in \cA$ and $0 \leq k \leq |\psi(ab)| - 2$, we
have
\[
  E\psi E(z)E(w) = E(\psi(y)w) = E(w\psi_k(y)) = E(w)E\psi_kE(z).
\]
Thus, there exists a word of length $k$, namely $w^\prime = E(w)$,
such that for any $u \in \cAstar$,
\[
  E\psi E(u)w^\prime = w^\prime E\psi_kE(u).
\]
\end{proof}
\begin{theorem} \label{07.07.03(2)}
Let $k, m \in \NN$ be such that $k = q_{m+1} - p$ with ~$2 \leq p
\leq q_{m+1} - q_{m} + 1$. Then
\[
  (\hat{\sigma}^m)_k(c_{1-\alpha}) = \prod_{j=m-1}^\infty
(\hat{\sigma}^{j+1})_{q_{m+1} - q_m - p}(a).
\] \vvqed
\end{theorem}

\section{Concluding Remarks}
Note that, by Corollary \ref{Cor:10.06.03(2)}, for any integer $m \geq 1$,
\[
  \sigma^m(b) = s_{m(n-1)}^{d_n - d_1}s_{m(n-1)-1},
\]
where $\sigma$ is the standard morphism that generates $c_\alpha$, for
$\alpha = [0;1+d_1,\ov{d_2,\cdots, d_n}]$ with $d_{n} \geq d_1 \geq 1$.
Now, by the periodicity of the continued fraction expansion of $\alpha$,
$d_i = d_{i+n-1}, ~\mbox{for all} ~i \geq 2$.
Hence, it is easily deduced that
\[
  d_n = d_{m(n-1) + 1}, ~~\mbox{for any integer $m \geq 1$},
\]
and one may write
\[
  \sigma^m(b) = s_{m(n-1)}^{d_{m(n-1)+1} - d_1}s_{m(n-1)-1}.
\]
Whence, if $d_1 = 1$, then for each $m \geq 1$,
\[
  v_{m(n-1)-1} = \begin{cases}
                a\sigma^{m}(b)b^{-1} ~~~\mbox{if $m$ is even}, \\
                b\sigma^{m}(b)a^{-1} ~~~\mbox{otherwise}.
                 \end{cases}  
\]
Accordingly, we do not have a `nice' expression for each $v_{k}$ ($k \in \NN$) 
in terms of
a power of $\sigma$ unless $n=2$, i.e. unless $\alpha = [0;2,\ov{r}]$ for
some $r \in \NN^{+}$. 
Therefore,
we cannot establish an extension of Theorem~\ref{T:12.06.03(1)} to the
case of $c_\alpha$ (nor $c_{1-\alpha}$) with $\alpha$ (as above) 
having $d_1 \geq 2$ and $n \geq 3$.

\section{Acknowledgements} The author would like
to thank Bob Clarke and Alison Wolff for their careful reading of the 
paper and useful suggestions.
The research was supported by the George Fraser Scholarship
of the University of Adelaide.


\end{document}